\documentclass[11pt,reqno]{amsart}
\usepackage[centertags]{amsmath}

\theoremstyle{plain} 

\newtheorem{lem}{Lemma}

\newtheorem{thm}{Theorem}

\theoremstyle{definition}
\newtheorem{defn}{Definition}
\newtheorem{example}{Example}

\theoremstyle{remark}
\newtheorem{rem}{Remark}  

\numberwithin{equation}{section}

 \DeclareMathOperator{\Ad}{Ad}
 \DeclareMathOperator{\opt}{opt}

\DeclareMathOperator{\rank}{rank}

\DeclareMathSymbol{\R}{\mathalpha}{AMSb}{"52}
\DeclareMathSymbol{\C}{\mathalpha}{AMSb}{"43}
\newcommand{\mbb}[1]{\mathbb{#1}}
\newcommand{\Z}{\mbb{Z}}

\newcommand{\K}{\mbb{K}}

\newcommand{\set}[1]{\left\{#1\right\}}

\newcommand{\ra}{\rightarrow}

\newcommand{\bd}{\begin{description}}
\newcommand{\ed}{\end{description}}
\newcommand{\beqr}{\begin{eqnarray}}
\newcommand{\eeqr}{\end{eqnarray}}
\newcommand{\beqt}{\begin{equation}}
\newcommand{\eeqt}{\end{equation}}

\begin{document}

\title[Invariants of Lie Algebras]{Invariants of a semi-direct sum of Lie algebras}
\author[J C Ndogmo]{J C  Ndogmo}

\address{
Department of Mathematics\\
 Univ. of Western cape\\
Private Bag X17\\
Bellville 7535\\
South Africa.}
\email{jndogmo@uwc.ac.za}

\begin{abstract}
We show that any semi-direct sum $L$ of Lie algebras with Levi
factor $S$ must be perfect if the representation associated with
it does not possess a copy of the trivial representation. As a
consequence, all invariant functions of $L$ must be Casimir
operators. When $S= \frak{sl}(2,\K),$ the number of invariants is
given for all possible dimensions of $L$. Replacing the
traditional method of solving the system of determining PDEs by
the equivalent problem of solving a system of total differential
equations, the invariants are found for all dimensions of the
radical up to five. An analysis of the results obtained is made,
and this lead to a theorem on invariants of Lie algebras depending
only on the elements of certain subalgebras.
\end{abstract}

\maketitle

\section{Introduction}

The subject of invariant functions of the coadjoint representation
of Lie algebras has given rise to many publications in recent
scientific literature \cite{abl, pavlow, per, ndogsix, campole}.
These are functions on the dual space $L^*$ of the Lie algebra $L$
that are invariant under the coadjoint action of the connected Lie
group $G$ generated by $L.$ \par

The invariants of  physical symmetry groups provide the quantum
numbers needed in the classification of elementary particles.
Thus, by making use of the eigenvalues of the Casimir operators of
the Poincar\'e group, Wigner achieved a classification of
particles according to their mass and spin \cite{wig}. One of the
most significant applications of the invariant functions in
physics is in the theory of dynamical symmetries, in which the
hamiltonian is written in terms of the Casimir operators of the
corresponding Lie symmetry group and its subalgebras
\cite{arim}.\par

For semisimple Lie algebras, these functions are well known,
following a paper by Racah \cite{rac}, which was mainly a
continuation of a work by Casimir \cite{casi} and other
physicists. For this class of Lie algebras, there always exists a
fundamental set of invariants that consists of homogeneous
polynomials, and the degree of transcendence of the associative
algebra generated by these functions over the base field $\K$ of
characteristic zero  equals the dimension of the Cartan
subalgebra.

The Levi decomposition theorem gives a preliminary classification
of Lie algebras into semisimple and solvable ones,  and a third
class consisting of  semi-direct sums of semisimple and solvable
Lie algebras. One fact that has emerged about these types of
semi-direct sums of Lie algebras is that all the invariants that
have been computed for them in the recent literature are
polynomials \cite{campole, pecina, per}. We consider in this paper
a semi-direct sum of Lie algebras of the form $L= S
\oplus_{\,\pi}\mathcal{R}$, where $S$ is semisimple and
$\mathcal{R}$ is the solvable radical, and where $\pi$ is a
representation of $S$ in $\mathcal{R}$ defining
the$[S,\mathcal{R}]$-type commutation relations. Unless otherwise
stated, we shall refer to a Lie algebra of this form simply as a
semi-direct sum of Lie algebras. The base field $\K$ is assumed to
be of characteristic zero. We show in section ~\ref{s: propinv}
that if $\pi$ does not possess a copy of the trivial
representation, then $L$ must be perfect; that is, equal to its
derived subalgebra. As a consequence of a result of \cite{abl},
this implies that all invariants of $L$ must be polynomials. Other
important consequences of this result are discussed. It is noted
in particular that all irreducible representations $\pi$ satisfy
the stated criteria.
\par

We move on to tackle in the next section the problem of the
explicit determination of the number of invariants of $L$ with
Levi factor $S= \frak{sl}(2, \K)$ and $\pi$ an irreducible
representation. Simple formulas giving this number are derived for
all possible dimensions of $\mathcal{R}.$ In section
~\ref{s:explicit}, we find it convenient to replace the most
common method that consists in solving a system of first order
PDEs for the determination of the invariants by the equivalent
problem of solving a system of total differential equations, as
described in a book by Forsyth \cite{fors}. This approach provides
a more tractable algorithm. The invariants are computed for all
dimensions of the radical up to 5, and an analysis of the
functions obtained is given. In particular, we show that when
$\dim \mathcal{R} > 3,$  all invariants depend only on the
elements of $\mathcal{R}.$ We show that this property always hold,
with a rank condition, for Lie algebras $L$ which are direct sums
$L_1 \oplus L_2$ of subspaces, and where $L_2$ is an abelian
subalgebra. Some families of Lie algebras having this property are
exhibited.\par
 The study of the invariants of Lie algebras reduces to
the case of indecomposable ones. We shall therefore assume, unless
otherwise stated, that the Lie algebra $L$ is indecomposable.

\section{Invariant functions and Casimir operators}

\subsection{Formal invariants}
Suppose that $L$ is the finite-dimensional Lie algebra of the
connected Lie group $G,$ and  denote  by $L^*$ the dual space of
$L.$ Let $\Ad \colon G \rightarrow GL(L) $ and $\Ad^* \colon G
\rightarrow GL(L^*)$ be the adjoint and the coadjoint
representations of $G$, respectively. Thus for $g\in G, f\in L^*$
and $x\in L$ we have $g\cdot f(x)=f( \Ad_{g^{-1}} (x)),$ where
$g\cdot f$ stands for $\Ad^*_g (f).$ Denote by $C^{\infty}(L^*)$
the space of all analytic functions on $L^*.$

\begin{defn}
A function $F\in C^{\infty}$ is called an invariant of the
coadjoint representation if $F(g\cdot f)=F(f),$ for all $g\in G$
and $f\in L^*.$
\end{defn}
Let $\{e_1, \dots, e_n\}$ be a basis of $L,$  and $\{f_1, \dots, f_n
\}$ the corresponding dual basis in $L^*$ given by $f_j(e_i)=
\delta_i^j.$ Let $\tilde{v}_i$ be the infinitesimal generator of the
coadjoint representation associated to the basis vector $e_i$ of
$L,$ and suppose that $x_1, \dots, x_n$ is a coordinate system in
$L^*$ corresponding to the dual basis. It is a well known fact
\cite{fom} that the infinitesimal generators are given by
\begin{equation}
\label{e: infge} \tilde{v}_i = - \sum_{j,\,k}  c_{ij}^k x_k
\frac{\partial}{\partial x_j} \qquad \text{(for $i=1,\dots, n$)}
\end{equation}
where the $c_{ij}^k$ are the structure constants of the
$n$-dimensional Lie algebra $L$ in the given basis. These vector
fields on $L^*$ form a Lie algebra that is homomorphic to $L.$ We
have more precisely, $[V_i, V_j] = \sum_k c_{ij}^k V_k.$ On the
other hand, by  considering the infinitesimal action of $\Ad^*$ on
$L^*,$ it is easy  to see \cite{fom} that a function $F\in
C^{\infty}(L^*)$ is an invariant of the coadjoint representation
if and only if it satisfies the system of partial differential
equations
\begin{equation}
\label{e: deteq} \tilde{v}_i \cdot F = 0, \qquad \text{(for
$i=1,\dots, n$)}
\end{equation}
The most common method for finding the invariants \cite{pavlow,
abl, ndogsix, campole, beltra} is by solving the system of linear
first order partial differential equations given by ~\eqref{e:
deteq}. In contrast with the polynomial functions in terms of
which the invariants of semisimple Lie algebras can always be
expressed, the solutions to equation ~\eqref{e: deteq} generally
involve various kind of functions, including rational and
logarithmic functions, as well as functions in arctan. They are
therefore usually referred to as formal invariants of $L.$ Because
any functional relations among the invariants yields another
invariant, they are determined by a maximal set of functionally
independent invariants. Such a set is called a fundamental set of
invariants. The number of invariants of $L$ usually refers to the
cardinality of this set. It is a well known fact \cite{beltra,
ndogsix} that the number $\mathcal{ N}$ of invariants of the Lie
algebra $L$ is given by
\begin{equation}
\label{e: number} \mathcal{ N }= \dim L - \rank (\mathcal{M}_L),
\end{equation}
where $\mathcal{M}_L = \left(  \sum_{k=1}^n c_{ij}^k x_k\right)$
is the matrix of the commutator table of $L.$

\subsection{Invariant polynomial functions.}
Let $V$ be a finite dimensional vector space. The symmetric
algebra $\frak{S}(V^*)$ is called the algebra of polynomial
functions on $V,$ and is often denoted \cite{hum} by
$\frak{P}(V)$. When a fixed basis $\{f_1, \dots, f_n\}$ of $V$ is
given, $\frak{P}(V)$ becomes identified with a the algebra of
polynomials in $n$ variables, $f_1, \dots, f_n.$ Thus the
coadjoint representation $\Ad^*_g$ of the Lie group $G$ acts
naturally on $\frak{P}(L),$ via  $(\Ad^*_g f)(x)=
f(\Ad_g^{-1}(x)),$ for $f \in \frak{P}(L),$ and we see that the
set denoted $\frak{P}(L)^I$ of all elements of $\frak{P}(L)$ fixed
by this action is precisely the algebra of polynomial invariants
of the coadjoint representation.\par

Let $\frak{A}(L)$ be the universal enveloping algebra of the Lie
algebra $L,$ and denote by $\frak{A}(L)^I$ its center. That is,
$\frak{A}(L)^I$ is the subset of elements commuting with all $x
\in \frak{A}(L),$ or equivalently, with all $x \in L.$ The
elements of $\frak{A}(L)^I$ are called Casimir operators, and they
commute with all elements of a representation. It follows by
Schur's Lemma that in any irreducible representation they are
represented by scalars. Since any automorphism $\sigma \colon L
\ra L$ extends uniquely to an automorphism of $\frak{A}(L),$ there
is an action of   $\Ad^* (G)$ on $\frak{A}(L),$ and this sends
$\frak{A}(L)^I$ onto itself. It can be shown (\cite{hum, dix})
that $\frak{A}(L)^I$ is precisely the set of all
$\Ad^*(G)$-invariants of $\frak{A}(L).$ Finally, it can be proved
\cite{abl,dix, rac} that the associative algebras $\frak{P}(L)^I$
and $\frak{A}(L)^I$ are algebraically isomorphic. That is, there
is a one-one correspondence between the polynomial invariants of
the coadjoint representation and Casimir operators. When $L$ is
semisimple, the space of all invariants of the coadjoint
representation is precisely $\frak{P}(L)^I.$  In this case, all
invariant functions of the coadjoint representation are
polynomials, and hence Casimir operators. We show in the next
section that this is also true for any nontrivial semi-direct sum
$L= S \oplus_{\,\pi}\mathcal{R}$ of Lie algebras, when the
representation $\pi$ does not possess a copy of the trivial
representation.

\section{Properties of the invariant functions}
\label{s: propinv}

Suppose that the finite dimensional Lie algebra $L$ over the field
$\K$ of characteristic zero is a semi-direct sum of the semisimple
Lie algebra $S$ and the solvable ideal $\mathcal{R}.$  That is, we
have the vector space direct sum
\beqt \label{e: dirsum}  L= S \overset{.}{+} \mathcal{R},  \eeqt
where $[S,S]=S, \; [\mathcal{R},\mathcal{R}]\subset \mathcal{R}$
and $[S,\mathcal{R}] \subset \mathcal{R}.$ Furthermore, we shall
always assume that this semi-direct sum is nontrivial, which means
that $[S,\mathcal{R}] \neq 0$ (otherwise $L$ is decomposable). The
given of the Lie algebra $L$ is equivalent to the given of the
semisimple Lie algebra $S$, the solvable algebra $\mathcal{R},$
and a representation $\pi$ of $S$ in $\mathcal{R},$ that defines
the $[S,\mathcal{R}]$-type commutation relations. Furthermore,
$\pi(z)$ must be a derivation of $\mathcal{R}$ for each $z \in S.$
Indeed, the Jacobi identity applied to $z \in S$ and $u, v \in
\mathcal{R}$ leads to
\beqt \label{e: derivation} \pi(z)\cdot[u,v]=[\pi(z)\cdot u, v] +
[u, \pi(z)\cdot v ], \eeqt
and justifies this assertion. Using equation ~\eqref{e:
derivation} and the solvability of $\mathcal{R},$ it is easy to
see that for any irreducible representation of $S$ in
$\mathcal{R},$  the radical $\mathcal{R}$ must be abelian.\par

 Let $ \mathcal{R}^S= \set{v \in \mathcal{R} \colon \pi(S)v = 0}.$
 In \cite{dix}, an element of the subspace $\mathcal{R}^S$ is called an
invariant of the $S$-module $\mathcal{R}.$ Denote also by
$\pi(S)\mathcal{R}$ the subspace of $\mathcal{R}$ generated by all
$\pi(s)\mathcal{R},$ where $s\in S.$ Because $S$ is semisimple,
$\mathcal{R}$ is the direct sum of $\mathcal{R}^S$ and
$\pi(S)\mathcal{R}.$

\begin{lem}
\label{lem: lperfect} Let $L= S\oplus_{\,\pi} \mathcal{R}$ be a
nontrivial semi-direct sum of the semisimple Lie algebra $S$ and
the solvable Lie algebra $\mathcal{R},$ and suppose that the
representation $\pi$ defines the $[S,\mathcal{R}]$-type
commutation relations.

\begin{description}
\itemsep = 1mm
\topsep = 2.5mm
\item[(a)] If $\pi$ does not possess a copy of the trivial
representation, then $L$ is perfect, and it has therefore a
fundamental set of invariants consisting of polynomials.
\item[(b)] The representation $\pi$ does not possess a copy of the trivial
representation if and only if $\pi(S)\mathcal{R}=\mathcal{R}.$

\end{description}
\end{lem}

\begin{proof}
We first notice that we have in this case $\pi(S)\mathcal{R}
=[S,\mathcal{R}].$ Now, for part (a), we see that $\pi$ does not
possess a copy of the trivial representation if and only if
$\mathcal{R}^S =0.$ It then follows from the remark preceding the
lemma that $[S,\mathcal{R}]=\mathcal{R},$ whence the equality
$[L,L]=L.$ The remaining part of the assertion is a consequence of
a result of \cite[Corollary 2]{abl} asserting precisely that any
perfect Lie algebra has a fundamental set of invariants consisting
of polynomial functions.
\par For part (b), the result is a
consequence of the fact that $\pi$ does not possess a copy of the
trivial representation if and only if $\mathcal{R}^S =0,$ and the
equality $ \mathcal{R}\\ = \pi(S)\mathcal{R} \oplus
\mathcal{R}^S.$
\end{proof}

\begin{thm}
\label{th: irred} Let $L= S\oplus_{\,\pi} \mathcal{R},$ with the
usual notation, and suppose that the representation $\pi$ is
irreducible. Then $L$ is a perfect Lie algebra, and has therefore
a fundamental set of invariants that consists of polynomial
functions.
\end{thm}

\begin{proof}
By part (a) of Lemma ~\ref{lem: lperfect}, we only need to show
that if $\pi$ is irreducible, then it does not have a copy of the
trivial representation. If $\pi$ is irreducible, then since the
semi-direct sum is nontrivial, the image space of a generic
element $\pi(z)$ of the representation is clearly a nonzero
invariant subspace, and its complementary subspace $W$ is the
largest subspace on which $\pi$ acts trivially. By irreducibility
$W=0,$ and thus $\pi$ does not possess a copy of the trivial
representation.
\end{proof}

\begin{rem}
\begin{enumerate}
\itemsep = 1mm
\topsep = 2.5mm

\item Not all semi-direct sums of Lie algebras are perfect. One
such example is given by the "optical Lie algebra" $\opt(2,1)$
\cite{zas}. It is a seven-dimensional subalgebra of the de Sitter
algebra $o(3,2)$ and has the form $L= S\oplus_{\,\pi}
\mathcal{R}$, where $S$ is generated by $\set{k_1, k_2, l_3},$ the
radical $\mathcal{R}$ is generated by $\set{w,m,q,c}$, and the
commutation relations are given by

\begin{alignat*}{5}
[w,m]&= - [k_1, m] = \textstyle{\frac{1}{2}}m \quad & [k_2,
q]&=[l_3,m]
=\textstyle{\frac{1}{2}}m \quad & [w,q] &= \textstyle{\frac{1}{2}}q \\
[k_1,q]&= [k_2, m]= \textstyle{\frac{1}{2}}q \quad & [w,c]&=
-[m,q] =c \quad &
-[l_3, m] &= \textstyle{\frac{1}{2}}q \\
 [k_1, k_2]&=-l_3 \quad & [k_1, l_3]&=-k_2 \quad & [k_2, l_3]& =
k_1.
\end{alignat*}

A simpler example is given by any Lie algebra of the form $L=
S\oplus_{\,\pi} \mathcal{R},$ where the radical $\mathcal{R}$ is
abelian and $\mathcal{R}^S \neq 0.$
\item For a Lie algebra of the form $L= S\oplus_{\,\pi} \mathcal{R}$, the
condition that $\pi$ does not possess a copy of the trivial
representation (or equivalently, that
$[S,\mathcal{R}]=\mathcal{R}$) is only a sufficient condition, but
not a necessary condition, for $L$ to be perfect. For example, the
derived subalgebra of $\opt(2,1)$ has $c$ as a central element,
and $\pi$ has therefore a copy of the trivial representation.
However, this subalgebra is perfect.

\item Lemma ~\ref{lem: lperfect} gives an interpretation of the
sufficient condition $[S,\mathcal{R}]=\mathcal{R}$ for $L$ to be
perfect in terms of the representation $\pi$ defining the
$[S,\mathcal{R}]$-type commutation relations.

\end{enumerate}
\end{rem}

By the already cited result of Racah \cite{rac}, any semisimple
Lie algebra has a fundamental set of invariants consisting of
polynomials. This result together with the Levi decomposition
Theorem and Lemma ~\ref{lem: lperfect}, shows that there are only
two types of Lie algebras that might not have a fundamental set
consisting of polynomial invariants. The first type consists of
semi-direct sums of Lie algebras of the form $L= S\oplus_{\,\pi}
\mathcal{R}$, where the representation $\pi$ possesses a copy of
the trivial representation (and is therefore not irreducible). The
second type consists of solvable non-nilpotent Lie algebras.
Indeed, by a result of \cite{bernat} the invariants of nilpotent
Lie algebras can all be chosen to be polynomials. However,
invariants in a fundamental set for solvable non-nilpotent Lie
algebras generally involve rational, logarithmic and other types
of functions \cite{ ndogsix,campori, pecina}. Although some
non-nilpotent solvable Lie algebras having a fundamental set
consisting of polynomials are given in \cite{pavlow}, no
characterization of such Lie algebras is available.\par For
semi-direct sums of Lie algebras, there is still no general result
concerning the number of their invariants, contrary to the case of
semisimple Lie algebras for which this number is known to be equal
to the rank of the algebra \cite{rac, hum}. We derive this number
for a particular Levi factor $S$ in the next section.

\section{The number of invariants of L}
\label{s: number} We suppose in this section that $S= \frak{sl}(2,
\K)$ has the standard basis $ x, y, h$ in which the commutation
relations are given by \beqt \label{e: sl2} [h, x]= 2 x, \qquad
[h, y]= - 2y, \qquad [x,y]=h. \eeqt
The following result is an immediate consequence of Theorem 7.2 of
\cite{hum} and Theorem 13.11 of \cite{sag}.

\begin{thm}
\label{th: repsl2} Let ${ S}= \frak{sl}(2,\K),$ where $\K$ is a
field of characteristic zero.
\begin{description}
\item[(a)] For any $m \in \Z^+,$ there exists a unique {\rm(}up to
isomorphism {\rm )} irreducible ${ S}$-module of highest weight
$m.$
\item[(b)] For each value $m\in \Z^+$ of the highest
weight, there exists a basis $\{v_0, \dots, v_m \}$ of the ${
S}$-irreducible module $V=V(m)$ in which the action of ${ S}$ is
given by
\begin{description}
\item[(b.1)] $h \cdot v_i = (m-2i)v_i$
\item[(b.2)] $y\cdot v_i = (i+1) v_{i+1}$
\item[(b.3)] $x\cdot v_{i} =  (m-i+1)v_{i-1}$
\end{description}

where $i\geq 0,$ $v_0$ is the maximal vector, and\\
 where $v_{-1}= v_{m+1}=0$ and $v_i= (1/i!) y^i \cdot v_0.$
\end{description}
\end{thm}

We shall also need the following result.

\begin{lem}
Let $M= \left( \begin{smallmatrix}
  A &B  \\
  C & 0
 \end{smallmatrix} \right)$ be a matrix partitioned into four blocks of
matrices $A,B,C$ and $0,$ where $0$ represents the zero matrix.
\begin{description}
\itemsep = 1mm
\topsep = 2.5mm

\item[(a)]If $A$ is a nonsingular square matrix of order $k,$ then
$M$ has rank $k$ if and only if $C A^{-1}B=0.$
\item[(b)] If $M$ is a nonsingular square matrix of order $2p$ so
that each of the block matrices is a square matrix of ~order $p,$
then $M^{-1}$ has the form $M^{-1} = \left( \begin{smallmatrix}
  0 & Y  \\
   Z & W
\end{smallmatrix} \right),$ where $Y,Z$ and $W$ are square
matrices of order ~$p.$
\end{description}

\end{lem}

\begin{proof}
Part (a) is an immediate consequence of a theorem from  a book by
Gantmacher \cite[P.47]{grant}. For part (b), we note that $M^{-1}
= \left( \begin{smallmatrix}
  X &Y  \\
  Z & W
 \end{smallmatrix} \right)$ implies $Z B =I$ and $X B=0,$ where $I$ is the
identity matrix. Thus $B$ is invertible and hence $X=0.$
\end{proof}

We now assume that, with the usual notation,  $ L = S \oplus_{\,
\pi} \mathcal{R}$  has finite dimension $n,$  and that $\dim
\mathcal{R} =d.$ Thus $L$ has dimension $n= d+3.$ We note that $d$
may assume any positive value, by part (a) of Theorem ~\ref{th:
repsl2}.

\begin{thm}
\label{rank} Let $L= S\oplus_\pi\! \mathcal{R},$ and set $r =
\rank (\mathcal{M}_L),$ where $\mathcal{M}_L$ is the matrix of the
commutator table of $L.$ Then
\begin{description}
\itemsep = 1mm
\topsep = 2.5mm
\item[(a)] $r=2,$ for $d=1.$
\item[(b)] $r=4,$ for $ d=2,3.$
\item[(c)] $r=6,$ for $d\geq 4.$
\end{description}
\end{thm}

\begin{proof}
For part (a), we note that when $d=1,$ $S$ acts trivially on the
one-dimensional module $\mathcal{R},$ and thus $\rank
(\mathcal{M}_L)= \rank (\mathcal{M}_S),$ where $\mathcal{M}_S$ is
the matrix of the commutator table of $S.$ It also follows from
equation ~\eqref{e: sl2} that $\mathcal{M}_S$ can be put into the
form
\begin{eqnarray*}
\mathcal{M}_S = \begin{pmatrix} A & B\\
C & 0
\end{pmatrix}, \quad \text{ with }\qquad  A= \begin{pmatrix} 0 & h \\
-h & 0 \end{pmatrix},
\end{eqnarray*}
and for some block matrices $B$ and $C,$ and we have $CA^{-1}B=0.$
The result then follows from Part (a) of
 the lemma.\par
 For part (b), one can always write $\mathcal{M}_L$ into the form
\begin{eqnarray*}
\mathcal{M}_L = \begin{pmatrix}
A_4 & B \\
C & 0 \end{pmatrix}, \qquad \text{ for } d=2,3,
\end{eqnarray*}
where  $A_4$ is a square non-singular matrix of order 4, and we
have $C A_4^{-1}B=0.$ Thus the result follows again from part (a)
of the lemma.\par Finally, we note that
\begin{eqnarray*}
\mathcal{M}_L = \begin{pmatrix}  M_6 & E \\ F & 0\end{pmatrix},
\quad \text{ for } \quad d\geq 4,
\end{eqnarray*}
where $M_6$ is a square matrix of order $6,$ of the form $\left(
\begin{smallmatrix}  A & B\\  C & 0 \end{smallmatrix} \right)$
and where $A, B$ and $C$ are $3 \times 3$ block matrices. Denoting
by $m=d-1$ the highest weight of the representation, we find that
$M_6$ has determinant
\begin{eqnarray*}
- (( -2+m) v_1 (( -1 + m) v_1^2 - 3 m v_0 v_2) + 3 m^2 v_0^2
v_3)^2
\end{eqnarray*}
and this is different from 0 since it has $9 m^4 v_0^4 v_3^2$ as a
term. Thus $M_6^{-1}$  has the form $M_6^{-1} = \left(
\begin{smallmatrix} 0_{3\times 3} & Y \\ Z & W
\end{smallmatrix}\right),$ by part (b) of the lemma. Noting now that the
sub-matrices $E$ and $F$ have their last three rows and last three
columns consisting of zeros, respectively, it follows that $F
M_6^{-1}E =0.$ By part (a) of the lemma again, we see that $\rank
(\mathcal{M}_L)= \rank (M_6)= 6 $
\end{proof}

\begin{thm}
\label{th: numbinv} Suppose that the radical $\mathcal{R}$ has
dimension $d.$ then the number ${\mathcal N}$ of invariants of $L=
S\oplus_\pi\! \mathcal{R}$ is given by
\begin{description}
\itemsep = 1mm
\topsep = 2.5mm
\item[(a)] ${\mathcal N}=2, \quad$ for $d=1$
\item[(b)] ${\mathcal N}=1,2, \quad$ for $d=2,3,$ respectively
\item[(c)] ${\mathcal N} = d-3, \quad $ for $d\geq 4 $
\end{description}
\end{thm}
\begin{proof}
It is a direct consequence of equation  ~\eqref{e: number} and
Theorem ~\ref{rank}.
\end{proof}

\section{Explicit determination}
\label{s:explicit}
\subsection{The method of total differential equations}
The invariants are usually determined by solving the system of
determining equations ~\eqref{e: deteq}. This is a system of
homogeneous linear first order partial differential equations, and
the method of characteristic is the most common for solving them.
These equations are suitable for the determination of the
invariants when the number of variables they involve is relatively
low. They are widely used for the determination of the invariants
\cite{pavlow,ndogsix, abl,campole}.\par

However, when the number of variables involved in the invariants
becomes relatively high, the equivalent adjoint system of total
differential equations becomes more appropriate for the
determination of the invariants. It provides a more efficient
algorithm involving a smaller number of change of variables and
substitutions. In particular, the number of equations in the
system corresponds to the number of invariants. Thus only one
total differential equation is to be solved when there is only one
invariant, no matter what the initial number of equations in
~\eqref{e: deteq}. We now derive the relationship between a system
of integral equations and the corresponding adjoint system of
total differential equations (see \cite{fors} for more details).

Let
\begin{equation}
\label{integraleq}
 \phi_j (u_1,\dots,u_p, x_1,\dots,x_q) = c_j
\end{equation}

be a system of p integral equations in $p+q$ variables, where $p$
and $q$ are positive integers. Assume furthermore that the
functions $\phi_j$ are functionally independent, for $j=1,\dots,
p.$ Then, without loss of generality, this amounts to assuming
that
\begin{equation}
\label{jacobian}
  \frac{\partial
(\phi_1,\dots,\phi_p)}{\partial (u_1,\dots,u_p)}  \neq 0
\end{equation}

Taking the differential in  ~\eqref{integraleq} yields
\begin{equation}
\label{differential}
 \sum_s \frac{\partial \phi}{ \partial u_s }
{\rm d} u_s + \sum_t \frac{\partial \phi}{\partial x_t} {\rm d}
x_t =0.
\end{equation}

By the implicit function theorem, condition ~\eqref{jacobian}
implies that one can solve the system ~\eqref{integraleq} for the
variables $u_j\quad  (j=1,\dots,p)$ in terms of the remaining $q$
variables. Thus we call the variables $u_j$ the dependent
variables and the remaining $q$ variables are the independent
variables of ~\eqref{integraleq}. In particular, ~\eqref{jacobian}
implies that one can uniquely solve for the differentials ${\rm
d}u_s$ in terms of the ${\rm d} x_t$ in ~\eqref{differential}.
This yields the following system of total differential equations.

\begin{equation}
\label{adjoint} {\rm d}u_s= \sum_{t=1}^q U_{s,t\,}{\rm d}x_t,
\qquad \textup{( for $s=1,\dots,p$).}
\end{equation}

The substitution of ~\eqref{adjoint} into the left hand side of
~\eqref{differential} gives rise to  a linear function of ${\rm d}
x_1, \dots,{\rm d}x_q.$ Because the variables $x_1,\dots,x_q$ are
independent, there cannot be any functional relation among them.
Thus the coefficients of the linear function in $ {\rm d} x_1,
\dots,{\rm d}x_q$ must all vanish. This leads to the following
system of $q$ linear first order partial differential equations,
called the Jacobian system of ~\eqref{integraleq}.

\begin{equation}
\label{jacobiansys} \Delta_t \phi \equiv \frac{\partial
\phi}{\partial x_t}+ \sum_{s=1}^p U_{s,\,t} \frac{\partial
\phi}{\partial u_s} = 0, \qquad \text{( for $t=1,\dots,q$).}
\end{equation}

It is easy to see \cite{fors} that the two equations
~\eqref{jacobiansys} and ~\eqref{adjoint} have exactly the same
integrals, and that a function $\Psi$ is a solution of
~\eqref{adjoint} or ~\eqref{jacobiansys} if and only if it is a
function of $\phi_1, \dots,\phi_p.$ It is also clear that both
equation ~\eqref{adjoint} and ~\eqref{jacobiansys} have exactly
$p$ functionally independent solutions. Equation ~\eqref{adjoint}
is called the adjoint system of ~\eqref{jacobiansys}. Conversely,
any equation of the form ~\eqref{jacobiansys} satisfies the
conditions of integrability if and only if the commutators of the
differential operators $\Delta_t \; \text{( for $t=1,\dots,q$)}$
all vanish \cite{peci16,fors}.

\subsection{Applications}
Suppose that the $n$-dimensional Lie algebra $L$ is generated by
$\{ X_1,\dots,X_n \}$ and that $\{x_1,\dots,x_n\}$ is a coordinate
system associated with this basis. We notice that under the
identification of the symmetric space $\frak{S}(L^*)$  with  $
\frak{S}(L)$, the coordinates system in a given basis of $L^*$ may
be replaced by the coordinates system in  the corresponding dual
basis in $L.$ Furthermore, it is customary to use the same
notation for basis vectors and corresponding coordinates systems
in the expression of the invariants.\par

By ~\eqref{e: deteq}, the determining equations ${\tilde X_i \cdot
F}=0, \text{ for $i=1,\dots n$ }$ are given by the system of
homogeneous linear fist order PDEs
\begin{equation} \label{deteq2} \sum_j [x_i, x_j]\cdot \frac{\partial
F}{\partial x_j} = 0, \qquad \text{ ( $i=1,\dots, n $),}
\end{equation}
where we have set  $[x_i, x_j]= \sum_k c_{ij}^{\,k} x_k.$
 Let $q= \rank (\mathcal{M}_L)= \rank ([x_i,x_j]).$
 Then we can solve
 ~\eqref{deteq2} for $q$ of the variables $\frac{\partial F}{\partial
 x_j}.$ This yields an equation of the form
 \begin{equation}
 \label{thejacobian}
\frac{\partial F}{\partial x_t} = - \sum_{s=q+1}^n U_{s,\,t\,}
\frac{\partial F}{\partial x_s}, \qquad (t = 1,\dots, q).
 \end{equation}
This equation represents the jacobian system for the integral
equations of ~\eqref{deteq2}.  Indeed, by equation ~\eqref{e:
number}, equation ~\eqref{deteq2} possesses exactly $n-q$
functionally independent invariants. It also determines the
coefficients $U_{s,\,t\,}$  of the corresponding adjoint system,
and thus the adjoint system itself, which is given by
\begin{equation}
\label{theadjoint} {\rm d} u_s = \sum_{t=1}^q U_{s,\,t\,} {\rm
d}x_t, \qquad \text{(for $s=1,\dots,n-q,$ and $u_s =x_{q+s}$)}.
\end{equation}
Consequently, the determining equations given by ~\eqref{deteq2}
and which represent a system of linear first order PDEs is
equivalent to ~\eqref{theadjoint}. Therefore, the latter system
can be used for the determination of the invariants. As already
noted, this system involves in many instances, and in particular
when the number of variables  in the equations is relatively high,
a more systematic algorithm than that needed for solving the
system of PDEs ~\eqref{deteq2} directly.

\subsection{Examples }
The method consisting of solving the system ~\eqref{theadjoint}
rather than the equivalent system of PDEs ~\eqref{deteq2} has been
used in \cite{pecina} for the determination of the invariants.
Namely, the invariants  of  certain solvable Lie algebras of
dimension six, and those of the Lie algebra $sa(n,\R)\; \text{(
for $n=2,3,4$)},$ where $sa(n, \R)$ is the semi-direct sum of
$\frak{sl}(n,\R)$ and the abelian Lie algebra of dimension $n.$ We
apply the same method here for the Lie algebra $L= S\oplus_\pi\!
\mathcal{R},$ where as in section ~\ref{s: propinv} this notation
represents the semi-direct sum of $S= \frak{sl}(2,\K),$ and the
radical $\mathcal{R},$ and where the commutation relations of the
$[S,\mathcal{R}]$-type are defined by the irreducible
representation $\pi.$\par

We shall suppose that the $n$-dimensional Lie algebra $L$ has a
basis of the form given in Theorem ~\ref{th: repsl2}. Thus
$\frak{sl}(2, \K)$ is generated by $\{x,y,h\}$ and the radical
$\mathcal{R}$ of dimension $d$ is generated by $\{v_0, v_1, \dots,
v_m\}$, where $m=d-1.$ It is clear from equation ~\eqref{deteq2}
that the determining equations are completely determined by the
commutator table of $L.$ In turn, this table is completely
determined by equation ~\eqref{e: sl2} and part (b) of Theorem
~\ref{th: repsl2}. However, to derive the adjoint system
~\eqref{theadjoint}, we also need the rank of the matrix of the
commutator table, and this is given by Theorem ~\ref{rank}. Once
the adjoint system is obtained, it can be solved using one of the
methods available for solving systems of total differential
equations, and Natani's method appears to be the most appropriate
in this case.
\par

Table 1 gives a list of invariants computed for $L=S\oplus_{\,\pi}
\mathcal{R}(m)$ where the dimension $d$ of the irreducible
$S$-module $\mathcal{R}(m)$ of maximal weight $m=d-1$ is up to
five, and where $L$ has the usual basis of the form $\{x,y,h,
v_0,v_1,\dots,v_m\}.$ We give here an example of the procedure of
calculations when $d=3.$ The matrix $\mathcal{M}_L$ of the
commutator table has the form

\begin{equation}
\label{ml1} \mathcal{M}_L = \begin{pmatrix}
0     &   h  & -2x   & 0   & 2 v_0 & v_1 \\
-h    &   0  & 2 y   & v_1 & 2 v_2  & 0   \\
2 x   & -2y  & 0     & 2v_0& 0     & -2 v_2\\
0     & - v_1& -2v_0 &  0  & 0     &0\\
-2v_0 & -2v_2&0      & 0   &0      &0\\
-v_1  & 0    & 2 v_2 & 0   &0      & 0
\end{pmatrix}
\end{equation}

The notation $[x_i,x_j]= \sum_k c_{ij}^{\,k} x_k$ used in equation
~\eqref{deteq2} shows that the determining equations can be read
off the commutator table. For instance, if we denote  by $F_u$
 the partial derivative $\partial F/\partial u$ of a function $F$
  with respect to the variable $u,$ the functions invariant
under the subgroup generated by $h$ are given by the third row of
~\eqref{ml1} as

\begin{eqnarray*}
{\tilde h}\cdot F \equiv 2 x F_x +
 - 2y F_y + 2 v_0 F_{v_0}
 -2 v_2 F_{v_2} =0
\end{eqnarray*}

However, given the matrix $\mathcal{M}_L$ of the commutator table,
and knowing that it has rank 4, we only need to solve the equation
$\mathcal{M}_L \cdot X=0,$ for $4$ of the components in the vector
\begin{eqnarray*}
X= \left( F_x, F_y, F_h, F_{v_0}, F_{v_1}, F_{v_2}\right)
\end{eqnarray*}
in terms of the remaining two others. This determines the Jacobian
system and the coefficients $U_{s,\,t\,},$ and hence the adjoint
system. More precisely, solving for four of the variables in terms
of the two variables $F_{v_1}$ and $F_{v_2}$ determines the
coefficients $U_{s,\,t\,},$ with $v_1$ and $v_2$ as dependent
variables. These transformations yield after simplification the
following system of two total differential equations.
\begin{align*}
\label{adjointd3}
 (hv_0 + v_1 x) {\rm d} v_1 &= -2 v_2 v_0 {\rm d}x
+ 2 v_0^2 {\rm
d}y \\
&- v_1 v_0 {\rm d}h + 2 (v_2 x + v_0 y) {\rm d}v_0\\
 (hv_0 + v_1 x){\rm d} v_2 &= - v_1 v_2 {\rm d}x+ v_0 v_1 {\rm
 d}y\\
& -\frac{v_1^2}{2}{\rm d}h + (v_1y -h v_2){\rm d}v_0
\end{align*}

Solving this system  gives the invariants
\begin{align*}
I_1= &\,h v_1 + 2 v_2 x - 2 v_0y\\
 I_2= &\,v_1^2 - h^2 v_1^2 - 4 v_0 v_2 - 4 h v_1 v_2 x - 4 v_2^2
 x^2\\
 &\,+ 4 h v_0 v_1 y + 8 v_0 v_2 x y - 4 v_0^2 y^2.
\end{align*}

After simplification, these invariants that we call again $I_1$
and $I_2$ are given by
\begin{align*}
I_1&= h v_1 + 2 v_2 x - 2 v_0y\\
 I_2&= v_1^2 - 4 v_0 v_2
\end{align*}

The computation of the invariants becomes more and more
impractical when the dimension of the radical is greater then 5.
This is partly due to the complications that arise with the
solving of the adjoint system, owing to the fast growing number of
terms and variables that appear in the invariants.

\begin{table}
\caption{\label{tab1} Invariants of the semi-direct sum $L=
S\oplus_\pi\! \mathcal{R}(m).$}
\begin{tabular}{ c c l}
\hline
  \rule[-0.07in]{0in}{6mm} Dimension of $\mathcal{R}(m)$& $\quad $
Algebra & $\quad$ Invariants \cr \hline
 &   &  \cr
1&$\frak{sl}(2,F)\, \oplus \langle v_0\rangle\quad $& $I_1= 4 xy +
h^2$\cr
 &                                          & $I_2= v_0$ \cr
 &   &  \cr
2&$\frak{sl}(2,F)\oplus_\pi\! \mathcal{R}(1) \quad$&  $I_1= v_1^2
x+ v_0 v_1 h - v_0^2 y$\cr
  &                               &                  \cr
3 & $\frak{sl}(2,F)\oplus_\pi\! \mathcal{R}(2)\quad$   & $I_1= h
v_1 + 2 v_2 x - 2 v_0 y $    \cr
 & & $I_2 = v_1^2 - 4 v_0 v_2$  \cr
 & & \cr
4 & $\frak{sl}(2,F)\oplus_\pi\! \mathcal{R}(3)\quad$ &  $I_1 = 2
v_0 v_2^3 - 9 v_0 v_1 v_2 v_3$\cr
 & & $\qquad + \frac{27}{2}\, v_0^2 v_3^2
   +\frac{1}{2}\, v_1^2 v_2^2 + 2 v_1^3 v_3$ \cr
   & & \cr
 5& $\frak{sl}(2,F)\oplus_\pi\! \mathcal{R}(4)\quad$  & $ I_1 =
  -12 v_0 v_4 + 3 v_1 v_3 - v_2^2 $  \cr
 &  &  $I_2 =
   27 v_0 v_3^2 - 9 v_1 v_2 v_3 + 27 v_1^2 v_4 $  \cr
 &  &  $\qquad-72 v_0 v_2 v_4 + 2 v_2^3$  \cr
\hline

\end{tabular}
\end{table}

\section{Properties of the invariants}

We notice that all the invariants of $L$  computed  for $\dim
\mathcal{R}(m)= 1,\dots,5$  and given in Table ~\ref{tab1} are all
polynomials as stipulated by Lemma ~\ref{lem: lperfect}. Moreover,
in all cases there is a fundamental system of invariants
consisting of homogeneous polynomials. On the other hand, for
$d\geq 4,$ we realize that the invariants  $F$ depend only on the
elements of the irreducible $S$-module $\mathcal{R}(m)$. That is,
they are of the form $F=F(v_0, v_1\dots,v_m ).$ The following
result generalizes this observation.

\begin{thm}
\label{th: invsl2} Let $L= S \oplus_{\,\pi} \mathcal{R},$ where
$\pi$ is an irreducible representation of\\ $S = \frak{sl}(2, \K)$
in the radical $\mathcal{R},$ and suppose that $ \dim \mathcal{R}
\geq 4.$
\begin{description}
\item[(a)] Every invariant of $L$ has the form
$F= F(v_0, \dots, v_m).$\\
 That is, $F$ does not depend on the
variables $x, y$ and $h$ associated with $\frak{sl}(2, \K).$
\item[(b)] The invariants of $L$ are all completely determined by
the $[S,\mathcal{R}]$-type commutation relations alone.
\end{description}
\end{thm}

\begin{proof}
Let $\tilde{V}_0, \tilde{V}_1$ and $\tilde{V}_2$ be the
infinitesimal generators corresponding to the basis elements $v_0,
v_1$ and $v_2.$ By the determining equations ~\eqref{e: deteq},
any invariant $F$ must satisfy the system of equations
$\tilde{V}_0 \cdot F = \tilde{V}_1 \cdot F = \tilde{V}_2 \cdot F
=0.$ Since $\mathcal{R}$ is abelian according to an earlier
remark, this system corresponds to a linear equation of the form
$A \cdot Z=0, $ where $Z$ is the vector $(F_x, F_y, F_h )$ and
where $A$ is the submatrix of the commutator matrix
$\mathcal{M}_L$ located between position $(4,1)$ and $(6,3).$ The
determinant of $A$ is $(-2+3m-m^2) v_1^3 + (- 6m + 3 m^2)v_0 v_1
v_2 - 3 m^2 v_0^2 v_3,$ which is clearly nonzero, and this proves
the first part of the theorem.\par For part (b), we notice that
because of the condition $F_x=F_y=F_h ~{=0} $ just proven in the
first part above, and the fact that $\mathcal{R}$ is abelian, the
infinitesimal generators corresponding to the basis elements of
$\mathcal{R}$ all reduce to zero. On the other hand, by equation
~\eqref{e: infge} the infinitesimal generator $\tilde{E}_i$
corresponding to a basis element $E_i \in \{x, y, h \}$ reduces
to
$$ \tilde{E}_i = \sum_{j=0}^m [e_i, v_j] \frac{\partial  }{\partial v_j}$$
where $e_i$ is the corresponding coordinate for $E_i.$ This proves
the assertion and  completes the proof of the theorem.
\end{proof}

Invariant functions that depend only on elements of a particular
subalgebra occur frequently in the study of Lie algebras. Theorem
~\ref{th: invsl2} represents only a particular case of a more
general framework in which such functions usually occur. Suppose
that the finite dimensional Lie algebra $\frak{M} = L_1
\overset{.}{+} L_2$ is a vector space direct sum of the subspace
$L_1$ and the abelian subalgebra $L_2.$  Let $\{X_1, \dots, X_t
\}$ be a basis of $L_2$ and extend it to a basis $\{X_1, \dots,
X_t, X_{t+1}, \dots, X_s\}$ of $\frak{M}.$  Denote by $B$ the
matrix
$\left( [x_i, x_j]\right)_{\begin{subarray}{l}i=1, \dots, t \\
j=t+1,\dots, s \end{subarray}},$ in which $ x_1, \dots, x_s$ is a
coordinate system in the given basis of $\frak{M}.$ We have the
following generalization of Theorem ~\ref{th: invsl2}.

\begin{thm}
\label{th: geninv} Suppose that $\dim L_2 \geq \dim L_1$ and that
the matrix $B$ is of maximal rank.

\begin{description}
\item[(a)] Every invariant of $\frak{M}$ is of the form $F= F(x_1, \dots,
x_t).$ That is, $F$ depends only on the elements of the abelian
subalgebra $L_2.$
\item[(b)] All invariants of $\frak{M}$ are completely determined
by the $[L_1, L_2]$-type commutation relations alone.
\end{description}
\end{thm}

\begin{proof}
Let $\tilde{X}_i$ be the infinitesimal generator of the coadjoint
action corresponding to $X_i.$  For every invariant $F,$ equation
~\eqref{e: deteq} implies in particular that $\tilde{X}_i \cdot ~F
= ~0 \;\\
 \text{ ( for $i= 1, \dots, t$ )}.$ The corresponding
system of PDE's can be written as a system of linear equations of
the form $B \cdot Z =0,$ where $Z$ is the vector $(F_{x_{t+1}},
\dots,
 F_{x_s}).$ The condition $\dim L_2 \geq \dim L_1$ means that $s-t \leq
 t$ and this ensures that when $B$ has maximal rank, it contains
 an invertible submatrix so that the linear system $B \cdot Z =0$
 implies $Z=0.$ The rest of the proof is similar to that given for
 Theorem ~\ref{th: invsl2}.
\end{proof}

\begin{example}
Take $\frak{M} = L_1 \overset{.}{+} L_2$ to be a solvable and
non-nilpotent Lie algebra having an abelian nilradical $L_2.$ By a
result of \cite{mubsix, ndogcan}, we have $\dim L_2 \geq \dim
\frak{M}/2.$ Furthermore, we showed in \cite{ndogcan} that the
corresponding matrix $B$ has maximal rank. It follows that
$\frak{M}$ always satisfies the hypothesis of the theorem and thus
the invariants of $\frak{M}$ depend only on the elements of the
nilradical,  and they are completely determined by the
$[L_1,L_2]$-type commutation relations alone (This is Theorem 3
and Corollary 1 of \cite{ndogcan}). The invariants of solvable
non-nilpotent Lie algebras of dimension six over $\R$ having
abelian nilradicals are computed in \cite{ndogsix}. None of them
has a fundamental set consisting of polynomials. Moreover, they
usually involve logarithms and functions in arctan.
\end{example}

\begin{rem}
Theorem ~\ref{th: geninv} greatly simplifies the determination of
the invariants, by reducing at least by half the number of
equations in the system of determining equations given by
~\eqref{e: deteq}, and by reducing by $\dim L_1$ the number of
independent variables in these equations.
\end{rem}

We showed that a semi-direct sum of Lie algebras is perfect when
the representation associated with it does not possess a copy of
the trivial representation. In such a case, semi-direct sums of
Lie algebras always have a fundamental set consisting of
polynomial invariants. We now make use of table ~\ref{tab1} to
show that despite these facts , and contrary to the case of
semisimple Lie algebras, the number of their invariants is not the
same as the dimension of the Cartan subalgebra. Indeed, when $\dim
\mathcal{R} =2,$ the Lie algebra $L= \frak{sl}(2, \K) \oplus_{\,
\pi} \mathcal{R}(2)$ has only one invariant. However, it is easy
to see that $L$ is split over the field $\K$ of characteristic
zero, and has Cartan subalgebra $ \K h \oplus \mathcal{R}(2),$
which has dimension 3.

\section{Conclusion}
In this paper, we considered a semi-direct sum of Lie algebras of
the form $L = S \oplus_{\, \pi}\mathcal{R},$ where  $\pi$ is a
representation of the semisimple Lie algebra $S$ in the radical
$\mathcal{R}$ that defines the $[S,\mathcal{R}]$-type commutation
relations. We showed that $L$ is perfect when $\pi$ does not
possess a copy of the trivial representation and that this
condition is equivalent to the requirement that $\pi(S)R=R.$ In
this case $L$ has a fundamental set of invariants consisting of
polynomials (Lemma ~\ref{lem: lperfect}). In particular, $L$ has
this property when $\pi$ is irreducible (Theorem ~\ref{th:
irred}). The number of invariants are given in Theorem ~\ref{th:
numbinv}, when $S = \frak{sl}(2, \K).$ Using a method of total
differential equations, we were able to determine the invariants
when the dimension of the radical is up to five (Table 1), and to
derive a theorem on certain properties of these invariants
(Theorem ~\ref{th: invsl2}), as well as a generalization of this
theorem (Theorem ~\ref{th: geninv}). Finally, we showed that,
although the Lie algebras we considered have a fundamental set of
invariants consisting of polynomials, the cardinality of this set
is generally not equal to the dimension of the Cartan subalgebra,
as it is in the well known case of semisimple Lie algebras.

\end{document}